\documentclass[10pt]{article}
\usepackage{times}
\usepackage{amsmath,amsfonts,amstext,amssymb,amsbsy,amsopn,amsthm,eucal,slashed,amstext,enumerate,accents}
\usepackage{color}
\usepackage{txfonts}
\usepackage{dsfont}
\usepackage{hyperref}
\usepackage{graphicx}   
\usepackage[colorinlistoftodos]{todonotes}

  \theoremstyle{plain}
 \newtheorem{theorem}{Theorem}

\newtheorem{proposition}{Proposition}

 \newtheorem{lemma}{Lemma}

 \newtheorem{remark}{Remark}

\theoremstyle{definition}

\newcommand{\pa}{\partial}

\newcommand{\Sn}{\mathbb{S}^n}
\newcommand{\Sp}{\mathbb{S}}

\newcommand{\sH}{\mathsf{H}}

\newcommand{\RR}{\mathbb{R}}

\newcommand{\ud}{\,\mathrm{d}}
\newcommand{\md}{\mathrm{d}}

\newcommand{\ve}{\varepsilon}

\allowdisplaybreaks

\title{
Critical nonlocal GJMS operator on odd dimensional spheres and mean field equation}
\author{Xuezhang Chen\thanks{$^\dag$X. Chen and J. Gan are partially supported by NSFC (No.11771204). Email: xuezhangchen@nju.edu.cn.},~ Jiaxue Gan\thanks{J. Gan's Email: jiaxuegan@smail.nju.edu.cn.}~~and Yalong Shi\thanks{Y. Shi is partially supported by NSFC (No.12371058). Email:shiyl@nju.edu.cn.}\\
		{\small $^{\ast}$$^{\dag}$$^{\ddag}$School of Mathematics \& IMS, Nanjing University, Nanjing 210093, P.R. China }}
\date{}

\begin{document}

\maketitle

\begin{abstract}
For a smooth axially symmetric function on $\mathbb{S}^n$ with $n\geq 1$, in odd dimensions we prove a neat  formula for the critical nonlocal GJMS operator $P_n$ similar to the one in even dimensions. With this formula, we give a simple proof of the Liouville theorem by Gui-Hu-Xie (J. Funct. Anal. 2022): Every smooth axially symmetric solution to
$$\frac{1}{n+1} P_n u = (n-1)!(e^{nu}-1) \quad \mathrm{on} \quad \Sn$$
is identically zero, without using any delicate property of Gegenbauer polynomials.

\medskip

{\bf MSC2020:} 35J30, 35J60, 35A22, 53C21.
\end{abstract}
\medskip

\section{Introduction}

On  $\Sn$, the GJMS operator in the critical dimension $n\geq 2$  is
\begin{equation}\label{def:critical-GJMS}
\begin{aligned}
P_n:=P_n^{\Sn}=\begin{cases}
\displaystyle  \prod_{j=0}^{(n-2)/2}\Big(-\Delta_{\Sn}+j(n-1-j)\Big) &\qquad  n \mathrm{~~even};\\
\displaystyle  \sqrt{-\Delta_{\Sn}+(\frac{n-1}{2})^2} \prod_{j=0}^{(n-3)/2}\Big(-\Delta_{\Sn}+j(n-1-j)\Big) &\qquad  n \mathrm{~~odd}.
\end{cases}
\end{aligned}
\end{equation}
As in Branson \cite[p.231]{Branson}, for any $n\in \mathbb{Z}_+$ we denote $B_n=\sqrt{-\Delta_{\Sn}+(\frac{n-1}{2})^2}$. Then the critical GJMS operator becomes
$$ P_n=  \prod_{j=0}^{n-1}\Big(B_n-\frac{n-1}{2}+j\Big).$$
This formula works also for $n=1$. See also Graham-Jenne-Mason-Sparling \cite{GJMS} and Beckner \cite{Beckner}.

For $n$ even, an astonishing formula of the GJMS operator acting on axially symmetric functions on $\Sn$ first appeared in Gui-Hu-Xie \cite{Gui-Hu-Xie}, which was established by a direct computation in lower dimensions. Chen-Shi \cite{Chen-Shi} gave a proof of this formula by using Gegenbauer polynomials.
\begin{proposition}\label{prop:Chen-Shi}
For $n$ even, suppose a smooth function $u$ on $\mathbb{S}^n$ depends only on $x:= x_{n+1}$. Then there holds
\begin{equation}\label{evenGJMS}
P_n u(x) = (-1)^{\frac{n}{2}}\left[ (1-x^2)^{\frac{n}{2}} u'\right]^{(n-1)}.
\end{equation}
\end{proposition}

For $n$ odd, it is natural to ask whether there exists an analogous formula of \eqref{evenGJMS} on $\Sn$. In this paper, we gave an affirmative answer. 

To continue, define the finite Hilbert transform on $(-1,1)$ by
\begin{equation}\label{Hilbert_transform}
\mathcal{H}[f](x) = \frac{1}{\pi} \mathrm{p.v.} \int_{-1}^1 \frac{f(t)}{x-t} \ud t, \quad -1 < x<1.
\end{equation}

For $n$ odd, we have the following
\begin{proposition}\label{prop:GJMS-odd-dim}
For $n=2m+1\geq 1$, suppose a smooth function $u$ on $\mathbb{S}^n$ depends only on $x:= x_{n+1}$. Then there holds
\begin{equation}\label{oddGJMS}
P_n u(x) = (-1)^m \left[ \mathcal{H}\big[ (1-t^2)^{\frac{n}{2}} u'(t)\big](x)\right]^{(2m)}.
\end{equation}
\end{proposition}
We will give a unified derivation of Propositions \ref{prop:Chen-Shi} and \ref{prop:GJMS-odd-dim} by induction on $n$, without using Gegenbauer polynomials.

\medskip

 As an application, we study the Livioulle theorem for the higher order mean field equation on spheres. 
 
 In \cite{Gui-Hu-Xie}, Gui-Hu-Xie proved the following Liouville theorem for axially symmetric solutions in dimension $n\geq 2$ to  
\begin{equation}\label{eq:mean-field-higher-order}
\frac{1}{n+1} P_n u = (n-1)!(e^{nu}-1) \quad \mathrm{on} \quad \Sn,
\end{equation} using delicate properties of Gegenbauer polynomials.
See also Gui-Wei \cite{Gui-Wei} for $n=2$. Dimension one case is exceptional, and a complete picture of  smooth solutions to the mean field equation 
\begin{equation}\label{MFE:circle}
P_1 u=\lambda(e^u-1) \qquad \mathrm{~~on~~} \Sp^1,~~ \lambda \in \RR_+
\end{equation}
lies in the combination of Osgood-Phillips-Sarnak \cite{OPS} and X. Wang \cite{Wang}; see Remark \ref{rk:MFE-circle}.

In \cite{Chen-Shi}, the authors combined the formula \eqref{evenGJMS} and integration by parts to give a simple proof of the above Liouville theorem for \eqref{eq:mean-field-higher-order} with $n$ even.
 
Now the formulae in Propositions \ref{prop:Chen-Shi} and \ref{prop:GJMS-odd-dim} allow us to  give a simple proof of \cite[Theorem 1.1]{Gui-Hu-Xie} without using Gegenbauer polynomials. 

\begin{theorem}\label{thm:oddGJMS}
For $n\geq 2$, every axially symmetric smooth solution of \eqref{eq:mean-field-higher-order} is identically  $0$.
\end{theorem}

\section{GJMS operators on axially symmetric functions}
Throughout the paper, the  round metric on $\Sn, n\geq 2$ is expressed as
$$
g_{\Sn} = (1-x^2)^{-1} \md x^2 + (1-x^2) g_{\mathbb{S}^{n-1}}.
$$
For a smooth axially symmetric  function $u$, without loss of generality, we assume $u=u(x)$. Set 
$$
L_n = -(1-x^2) \pa_x^2 + nx \pa_x,
$$
then $-\Delta_{\Sn} u(x) = L_n u(x)$, and
$$
B_n =  \sqrt{L_n+(\frac{n-1}{2})^2}. 
$$
\begin{lemma}
For a smooth axially symmetric function $u=u(x)$ on $\Sn$, we have
\begin{equation}\label{commuteD&L}
    \pa_x L_n u = (L_{n+2}+n)\pa_x u.
\end{equation}
Consequently,
\begin{equation}\label{commuteD&A}
    \pa_x B_n^2 u = B_{n+2}^2\pa_x u.
\end{equation}
\end{lemma}
\begin{proof}
Direct differentiation gives
\begin{align*}
\pa_x L_n u &= \pa_x(-(1-x^2)u''+nxu')\\
&= -(1-x^2)u'''+2xu''+nxu''+nu'\\
&= (L_{n+2}+n) \pa_x u ,
\end{align*}
thereby
\begin{align*}
\pa_x B_n^2 u  &= \pa_x \left( L_n+(\frac{n-1}{2})^2 \right) u\\
&= \left( L_{n+2} +n+(\frac{n-1}{2})^2 \right) \pa_x u\\
&= \left( L_{n+2}+(\frac{n+1}{2})^2 \right) u \\
&= B_{n+2}^2 \pa_x u.
\end{align*}
\end{proof}

\begin{lemma}\label{lem:commuteD&B}
For a smooth axially symmetric  function $u=u(x)$, we have
\begin{equation*}
    \pa_x B_n u = B_{n+2} \pa_x u.
\end{equation*}
\end{lemma}
\begin{proof}
Write $A_n  = B_n^2$. Then by $(\ref{commuteD&A})$ we have
$$
\pa_x A_n u  = A_{n+2}\pa_x u.
$$
Fix $t>0$. Let $f = (t+ A_n)^{-1}u$, then
$$
\pa_x u = \pa_x \Big((t+ A_n)f \Big) = (t+ A_{n+2})\pa_x f.
$$
Therefore,
$$
(t + A_{n+2})^{-1}\pa_x u = \pa_x (t+ A_n)^{-1}u.
$$
Then
\begin{equation}\label{An}
\begin{aligned}
A_{n+2}(t + A_{n+2})^{-1}\pa_x u &= \Big(I - t (t+ A_{n+2})^{-1} \Big) \pa_x u \\
&= \pa_x u - t \pa_x (t + A_n)^{-1}u \\
&= \pa_x A_n (t + A_n)^{-1}u.
\end{aligned}
\end{equation}
Using the integral
$$
\sqrt{\lambda} = \frac{1}{\pi} \int_0^\infty t^{-\frac{1}{2}} \frac{\lambda}{t+\lambda} \ud t,  \quad \lambda \in \RR_+,
$$
we have
$$
B_n u = A_n^{\frac{1}{2}} u = \frac{1}{\pi} \int_0^\infty t^{-\frac{1}{2}} A_n (t+A_n)^{-1} u \ud t.
$$
So,
\begin{align*}
\pa_x B_n u &= \frac{1}{\pi} \int_0^\infty t^{-\frac{1}{2}} \pa_x A_n (t+A_n)^{-1} u \ud t \\
&= \frac{1}{\pi}  \int_0^\infty t^{-\frac{1}{2}} A_{n+2}(t + A_{n+2})^{-1}\pa_x u \ud t \\
&= B_{n+2} \pa_x u.
\end{align*}
Here, we use $(\ref{An})$ in the second equality.
\end{proof}

\begin{proposition}
For every integer $n \geq 1$, we have 
\begin{equation}\label{recurrenceP_n}
P_{n+2} \pa_x = L_{n+2} \pa_x P_n = \pa_x P_n(L_n -n).
\end{equation}
\end{proposition}
\begin{proof}
Introduce the scalar polynomial
$$
p_n (z) = \prod_{j=0}^{n-1} (z + j-\frac{n-1}{2}), \quad n \in \mathbb{Z_+}.
$$
The polynomials themselves satisfy the elementary recurrence
$$
p_{n+2}(z) = \Big(z^2 - \frac{(n+1)^2}{4}\Big)p_n(z).
$$
Then for all $n \in \mathbb{Z}_+$, the GJMS operator in \eqref{def:critical-GJMS} can be written in a unified way as
$$
P_n = p_n (B_n).
$$
Lemma \ref{lem:commuteD&B} gives $\pa_x p_n(B_n) = p_n (B_{n+2}) \pa_x$. Consequently,
\begin{align*}
P_{n+2} \pa_x &= p_{n+2}(B_{n+2}) \pa_x\\
&= \Big(B_{n+2}^2 - \frac{(n+1)^2}{4}\Big)p_n(B_{n+2}) \pa_x\\
&= L_{n+2}\pa_x p_n({B_n}) \\
&= L_{n+2}\pa_x P_n.
\end{align*}
Notice that $P_n$ commutes with $L_n$. Hence $(\ref{commuteD&L})$ gives
$$
P_{n+2} \pa_x=L_{n+2}\pa_x P_n = \pa_x (L_n -n)P_n = \pa_x P_n(L_n -n).
$$
\end{proof}

We now start with the GJMS operator $P_1$ on the unit circle $\mathbb{S}^1$: Under the polar coordinates $(r,\theta)$ such that $x_1 = \sin \theta, x=x_2 = \cos \theta$,  $-\Delta_{\mathbb{S}^1} = - \pa_\theta^2$, then $P_1 = (- \pa_\theta^2)^{\frac{1}{2}} = \vert \pa_\theta \vert$.
This can be interpreted as the Dirichlet-to-Neumann map in the unit disk $\mathbb{D}_1$ (cf. Osgood-Phillips-Sarnak \cite{OPS}): Let $u$ be a smooth solution to
\begin{align*}
\begin{cases}
\Delta u=0 \qquad\quad \mathrm{~~in~~}\mathbb{D}_1,\\
u(1,\theta)=u(\theta).
\end{cases}
\end{align*}
If we expand 
$$u(\theta)=\sum_{k=-\infty}^\infty u_k e^{ik\theta},$$
then we have
$$u(r,\theta)=\sum_{k=-\infty}^\infty u_k r^{|k|} e^{ik\theta}$$
and define
$$P_1 u:=\pa_r u |_{r=1}=\sum_{k=-\infty}^\infty |k| u_k e^{ik\theta}.$$

Let $\mathcal{H}_{\mathrm{circ}}$ be the periodic Hilbert transform with Fourier multiplier $-i\mathrm{sgn}(k)$ on $e^{ik\theta}$.
Using the fact that for $x \notin 2\pi \mathbb{Z}$
$$
\cot \frac{x}{2} = 2 \sum_{k=0}^\infty \sin(kx) = -i \sum_{k=-\infty}^\infty \mathrm{sgn}(k) e^{ikx},
$$
we have
$$
\mathcal{H}_{\mathrm{circ}} f(\theta) = \frac{1}{2\pi} \mathrm{p.v.} (f*\cot \frac{x}{2}) =  \frac{1}{2\pi} \mathrm{p.v.} \int_{-\pi}^\pi f(\varphi) \cot \frac{\theta - \varphi}{2} \ud \varphi.
$$
Then we have
$$
P_1 u = \mathcal{H}_{\mathrm{circ}}\pa_\theta u.
$$
Since $\pa_\theta u = -\sin\theta u'(\cos \theta)$, we compute
\begin{align*}
P_1 u &= - \frac{1}{2\pi} \mathrm{p.v.} \int_{-\pi}^\pi \sin\varphi u'(\cos \varphi) \cot \frac{\theta - \varphi}{2} \ud \varphi\\
&= \frac{1}{2\pi} \mathrm{p.v.} \int_0^\pi \sin\varphi u'(\cos \varphi)\Big( \cot \frac{\theta + \varphi}{2} - \cot \frac{\theta - \varphi}{2} \Big) \ud \varphi \\
&= \frac{1}{\pi} \mathrm{p.v.} \int_0^\pi \frac{\sin^2\varphi u'(\cos \varphi)}{\cos\theta - \cos \varphi}\ud \varphi.
\end{align*}
Substitution $t = \cos \varphi$ gives
\begin{equation}\label{P_1}
P_1 u(x) =  \frac{1}{\pi} \mathrm{p.v.} \int_{-1}^1 \frac{\sqrt{1-t^2}u'(t)}{x-t}\ud t = \mathcal{H}[(1-t^2)^{\frac{1}{2}}u'(t)](x).
\end{equation}

Now we are ready to give a unified proof for $(\ref{evenGJMS})$ and $(\ref{oddGJMS})$.
\begin{proof}[Proof of Propositions \ref{prop:Chen-Shi} and \ref{prop:GJMS-odd-dim}.]
Define 
\begin{equation*}
\varepsilon_n =(-1)^{\lfloor \frac{n}{2}\rfloor}
\quad\mathrm{and~~}\quad
T_n =
\begin{cases}
\mathrm{Id}, & n \mathrm{~~even},\\
\mathcal H, & n \mathrm{~~odd}.
\end{cases}
\end{equation*}
Then our goal is to prove
\begin{equation}\label{P_n}
    P_n u = \varepsilon_n \pa_x^{n-1}T_n [(1-x^2)^{\frac{n}{2}}u'(x)].
\end{equation}

We prove it by induction. It's true for $n=1$ by \eqref{P_1} and for $n=2$ by a direct computation. Now assume $(\ref{P_n})$ for some $n \geq 1$. Take a smooth primitive $\mathsf{U}$ of $u$ such that $\mathsf{U}' = u$. By the recurrence $(\ref{recurrenceP_n})$ and \eqref{commuteD&L} we obtain
\begin{align*}
P_{n+2} u &= P_{n+2} \pa_x \mathsf{U} \\
&= \pa_x P_n (L_n - n) \mathsf{U} \\
&= \varepsilon_n \pa_x^n T_n [(1-x^2)^{\frac{n}{2}}\pa_x(L_n -n) \mathsf{U}]\\
&= \varepsilon_n \pa_x^n T_n [(1-x^2)^{\frac{n}{2}}L_{n+2}u] \\
&= -\varepsilon_n \pa_x^n T_n \pa_x[(1-x^2)^{\frac{n+2}{2}}u']\\
&= \varepsilon_{n+2} \pa_x^{n+1} T_{n+2} [(1-x^2)^{\frac{n+2}{2}}u'].
\end{align*}
Therefore, starting from $n=1$ proves all odd dimensions and starting from $n=2$ proves all even dimensions.
\end{proof}

\section{Mean field equations on spheres}

To prove Theorem \ref{thm:oddGJMS}, we first recall basic properties of Hilbert transforms and Fourier transforms.

We use the following convention for Fourier transform
$$
\hat{f}(\xi):=\mathcal{F}[f](\xi) = \frac{1}{\sqrt{2\pi}} \int_\RR e^{-ix\xi} f(x) \ud x.
$$

Distinguish the real line Hilbert transform by
\begin{equation*}
\mathsf{H}[F](x) =  \frac{1}{\pi} \mathrm{p.v.} \int_\mathbb{R} \frac{F(t)}{x-t} \ud t:=\mathrm{p.v.} K\ast F(x),
\end{equation*}
where $K(x)=\frac{1}{\pi x}$ for $x \in \RR\setminus \{0\} $.
In particular, for a function supported in $[-1,1]$, its restriction to $(-1,1)$ is exactly \eqref{Hilbert_transform}.

For brevity, let $T:= (-1)^m \frac{\ud^{2m}}{\ud x^{2m}}\circ\sH$ when $n=2m+1$, and $T:=(-1)^m \frac{\ud^{2m-1}}{\ud x^{2m-1}}$ when $n=2m$.
Then we have the following:

\begin{lemma}
We have $\widehat{\sH[f]}= - i\ \mathrm{sgn}(\xi) \hat{f}$, and
$$
\mathcal{F}[Tf] (\xi)= -i\xi\vert \xi \vert^{n-2} \hat{f}=:t(\xi)\hat{f}.
$$
\end{lemma}
\begin{proof}
The property of the Fourier transform gives $\widehat{\sH[f]}=\widehat{K \ast f}=\sqrt{2\pi} \hat K \hat f$. To compute $\hat K$, approximate $K$ by
$$K_{\ve}(x)=\frac{1}{\pi x}\mathbf{1}_{\{|x|> \ve\}}, \quad \ve \in \RR_+.$$
Then we have
\begin{align*}
\widehat{K_\ve}(\xi)=&\frac{1}{\sqrt{2\pi}}\frac{1}{\pi}\left(\int_{\ve}^\infty+\int_{-\infty}^{-\ve}e^{-ix \xi} \frac{ \ud x}{x}\right)=-\frac{1}{\sqrt{2\pi}}\frac{2i}{\pi}\int_\ve^\infty \frac{\sin x\xi}{x} \ud x,
\end{align*}
thereby,
$$\hat K(\xi)=\lim_{\ve \searrow 0}  \widehat{K_\ve}(\xi)=-\frac{i}{\sqrt{2\pi}} \mathrm{sgn}(\xi).$$
Hence, we obtain $\widehat{\sH[f]}=- i\ \mathrm{sgn}(\xi) \hat{f}$.

Notice that when $n=2m+1$,
$$\mathcal{F}[Tf] =(-1)^m(i\xi)^{2m} \widehat{\sH [f]}(\xi)=-i\xi\vert \xi \vert^{2m-1} \hat{f};$$
when $n=2m$,
$$\mathcal{F}[Tf] =(-1)^m(i\xi)^{2m-1} \hat{f}(\xi)=-i\xi\vert \xi \vert^{2m-2} \hat{f}.$$
This completes the proof.
\end{proof}

 A basic fact will be frequently used: Assume that $f \in C^1([-1,1])\cap C^\infty((-1,1))$, $f'' \in L^1 (-1,1)$, and
$$
f(\pm 1) = f'(\pm 1) = 0.
$$
Since for a Schwartz function  $h$, 
$$
(\sH [h])' = \sH[h'],
$$
Then applying $\sH$ to the zero-extended function $f$ on $\RR$ gives
$$
(\mathcal H [f])' = \mathcal{H}[f'].
$$

In the following, we always set $m=\lfloor \frac{n}{2} \rfloor$ for $n\geq 2$.
\begin{lemma}\label{lem:key}
Fix $n\geq 2$. Assume that $f \in H^m(\RR)$ and there exist $f_j \in C_c^\infty(-1,1)$ such that
$$
f_j \to f \quad \mathrm{~~in~~}  H^m(\RR) \quad \mathrm{~~as~~} j \to \infty.
$$
Then 
\begin{equation}\label{energy}
\langle Tf,xf \rangle_{H^{-m},H^m} = \frac{n-1}{2} \int_\RR |\xi|^{n-2}\hat |f(\xi)|^2 \ud \xi.
\end{equation}

\end{lemma}
\begin{proof}
First assume that $f \in C_c^\infty (-1,1)$. Since the multiplier $t(\xi)$ is purely imaginary and odd, the operator $T: H^m(\RR) \to H^{-m}(\RR)$ is skew-adjoint. Using the property of the Fourier transform that $\widehat{xf} = i\pa_\xi \hat{f}$, we obtain
\begin{equation}\label{Lie}
\begin{aligned}
\mathcal{F}[[T,x]f] =& \mathcal{F}[T(xf)-xTf]\\
=& t(\xi) \widehat{xf} - i\pa_\xi \widehat{Tf}\\
=& t(\xi)i\pa_\xi \hat{f} - i\pa_\xi(t(\xi) \hat{f}) \\
=& -i t'(\xi) \hat{f}= - (n-1)\vert \xi \vert^{n-2} \hat{f}.
\end{aligned}
\end{equation}

By Plancherel theorem and $(\ref{Lie})$ we obtain
\begin{align*}
\langle Tf,xf \rangle_{H^{-m},H^m} 
=& \int_\RR (Tf)\cdot(xf) \ud x\\
=& -\int_\RR f \cdot T(xf) \ud x\\
=& -\int_\RR f\cdot xTf \ud x -\int_\RR f [T,x] f \ud x\\
=&-\int_\RR xf  Tf \ud x-\int_\RR  \mathcal{F}[[T,x]f] \bar{\hat f} \ud \xi\\
=& - \left\langle Tf,xf \right\rangle_{H^{-m},H^m} + (n-1) \int_\RR |\xi|^{n-2}|\hat f(\xi)|^2 \ud \xi.
\end{align*}
This proves $(\ref{energy})$ for smooth compactly supported functions. 

For general $f \in H^m(\RR)$, we have
\begin{align*}
\lVert Tf\rVert_{H^{-m}}^2&= \int_\RR (1 + \vert\xi\vert^2)^{-m} \vert t(\xi) \vert^2 \vert \hat{f}(\xi)\vert^2 \ud \xi \\
&= \int_\RR (1 + \vert\xi\vert^2)^{-m} \vert \xi \vert^{4m} \vert \hat{f}(\xi)\vert^2 \ud \xi \\
&\leq \int_\RR (1 + \vert\xi\vert^2)^m \vert\hat{f}(\xi)\vert^2 \ud \xi \\
&= \lVert f(\xi)\rVert_{H^m}^2.
\end{align*}
Then $T: H^m(\RR) \to H^{-m}(\RR)$ is continuous. Moreover, since $f_j \to f$ in $H^m(\RR)$ and $\frac{n-2}{2} < m$, we obtain that $xf_j \to xf$ in $H^m(\RR)$ and $(1+\vert \xi )\vert^{\frac{n-2}{2}} \widehat{f_j} \to (1+\vert \xi\vert)^{\frac{n-2}{2}} \hat f$ in $L^2(\RR)$. Therefore, passing to a limit we prove \eqref{energy}.
\end{proof}

For a smooth axially symmetric function $u$, define its weighted derivative, extended by zero, 
\begin{equation}\label{zeroextension}
F(x)=\left\{  
     \begin{array}{lr}  
     (1-x^2)^{\frac{n}{2}}u'(x), & \quad\vert x \vert <1, \\
     0, &\quad  \vert x \vert \geq 1.
     \end{array}  
\right.  
\end{equation} 
Then it is not hard to show that $F \in H^{m}(\RR)$ and there exist $F_j \in C_c^\infty(-1,1)$ such that
$$
F_j \to F \quad \text{in} \quad H^m(\RR).
$$
Set 
$$
 U(x) = (1-x^2)u'(x).
$$

\begin{lemma}\label{lem:energy}
For $n\geq 2$, every  smooth axially symmetric function $u=u(x)$ satisfies
\begin{equation*}
\int_{\Sn} xU P_n u \ud \mu_{\Sn} = \frac{n-1}{2} |\Sp^{n-1}|\int_\RR |\xi|^{n-2}|\hat F|^2 \ud \xi.
\end{equation*}
\end{lemma}
\begin{proof}
By Propositions \ref{prop:Chen-Shi} and \ref{prop:GJMS-odd-dim} we obtain $P_n u=T F$ . Then by Lemma \ref{lem:key} we have
\begin{align*}
\int_{\Sn} xU P_n u \ud \mu_{\Sn}=&|\Sp^{n-1}|\int_{-1}^1 x(1-x^2)^{\frac{n-2}{2}}U\cdot TF \ud x\\
=&|\Sp^{n-1}|\langle TF,xF \rangle_{H^{-m},H^m}= \frac{n-1}{2} |\Sp^{n-1}|\int_\RR |\xi|^{n-2}|\hat F|^2 \ud \xi.
\end{align*}
\end{proof}

\begin{proof}[Proof of Theorem \ref{thm:oddGJMS}.]
Let
$F_2 = 1 - (n+1)x^2$ be a spherical harmonic on $\Sn$ of degree $2$.
Then
\begin{equation}\label{F2}
\Delta F_2 = -2(n+1)F_2, \quad P_n F_2 = (n+1)! F_2.
\end{equation}
Define
$$
J= \int_{\Sn} xU \ud \mu_{\Sn}.
$$
Then integration by parts yields
\begin{equation}
\begin{aligned}
J &= |\Sp^{n-1}| \int_{-1}^1 x(1-x^2)^{\frac{n}{2}} u'(x) \ud x \\
&= -|\Sp^{n-1}| \int_{-1}^1 u(x) F_2(x) (1-x^2)^{\frac{n-2}{2}} \ud x \\
&= - \int_{\Sn} uF_2 \ud \mu_{\Sn}.
\end{aligned}
\end{equation}

By  $(\ref{F2})$  we obtain
$$
\frac{1}{n+1} \int_{\Sn} F_2 P_n u \ud \mu_{\Sn}=\frac{1}{n+1} \int_{\Sn} u P_n F_2   \ud \mu_{\Sn} = -n! J.
$$
On the other hand, multiplying both sides of \eqref{eq:mean-field-higher-order} by $F_2$ and integrating over $\Sn$ yield
\begin{align*}
& \frac{1}{n+1} \int_{\Sn} F_2 P_n u \ud \mu_{\Sn}\\
=& (n-1)! \int_{\Sn} (e^{nu} -1) F_2  \ud \mu_{\Sn} \\
=& -\frac{(n-1)!}{2(n+1)} \int_{\Sn} (e^{nu} -1) \Delta F_2  \ud \mu_{\Sn} \\
=& \frac{n!}{2(n+1)} \int_{\Sn} e^{nu} \left\langle \nabla F_2, \nabla u \right\rangle \ud \mu_{\Sn} \\
=& -n! \int_{\Sn} e^{nu}x U \ud \mu_{\Sn} \\
=& - n! \int_{\Sn}  xU\left( \frac{1}{(n+1)(n-1)!}P_n u +1\right) \ud x \\
=& -\frac{n}{n+1} \int_{\Sn} xU P_n u \ud \mu_{\Sn} - n! J.
\end{align*}
Hence, putting these facts together yields
$$
\int_{\Sn} xU P_n u \ud \mu_{\Sn} = 0.
$$
Then Lemma \ref{lem:energy} gives
$$
0 =  \int_{\RR} \vert \xi \vert^{n-2} \vert \hat{F}(\xi) \vert^2 \ud \xi.
$$
It follows that $\hat{F} = 0$ a.e. in $\RR$, and then $F = 0$ a.e. in $(-1,1)$. Since $F$ is continuous in $[-1,1]$, then $F=0$ everywhere and thus $u'=0$. So $u$ is a constant.
Going back to equation \eqref{eq:mean-field-higher-order}, this constant equals $0$.
\end{proof}

\begin{remark}\label{rk:MFE-circle}
Consider the mean field equation \eqref{MFE:circle} on $\Sp^1$,when $\lambda \notin \mathbb{Z}$, the smooth solution to \eqref{MFE:circle} is identically zero by \cite[Lemma 2.3]{OPS}; when $\lambda=N \in \mathbb{Z}_+$, the smooth solutions to \eqref{MFE:circle} are 
$$ u(z)=\log \frac{1-|\xi|^2}{|1-\xi z^N|^2}, \quad z \in \mathbb{D}_1 \subset \mathbb{C}$$
 for some $\xi \in \mathbb{C}$ and $|\xi|<1$, by \cite[Theorem 3]{Wang}. The distinction is clear in the proof of Lemma \ref{lem:key}, where we use the derivative of $t(\xi)$. When $n=1$, $t(\xi)$ is the sign function and is not differentiable. There will be a $\delta$-function instead.
\end{remark}

\appendix
\section{Alternative proof of Proposition \ref{prop:Chen-Shi}}\label{Append}

The purpose is to give a direct induction proof, also independent of Gegenbauer polynomials.

\begin{proof}[Proof of Proposition \ref{prop:Chen-Shi}.] 
For $n$ even, let $X = \partial_x \left( (1-x^2)^{\frac{n}{2}} \partial_x \right)$ and 
$$e_l = \Pi_{j=0}^l (\Delta_{\Sn} - j(n-j-1)),\quad 0\leq l \leq \frac{n}{2}.$$  Then for each $l$, there holds 
 \begin{equation}\label{eqn:P_n-induction}
 e_l u= \partial_x^l \left[ (1-x^2)^{ l+1-\frac{n}{2}} \partial_x^l Xu \right].
 \end{equation}

We prove by induction on $l$. For $l=0$, it is true. Assume that \eqref{eqn:P_n-induction} holds for $l$. Then we have 
\begin{align*}
	&\partial_x^{l+1} \left[ (1-x^2)^{ l+2-\frac{n}{2}} \partial_x^{l+1} Xu \right] \\
	=& \partial_x^{l+1} \big[ (1-x^2) \partial_x  \left( (1-x^2)^{l+1 -\frac{n}{2}} \partial_x^l Xu \right)
	+2 (l+1-\frac{n}{2}) x (1-x^2)^{l+1-\frac{n}{2}} \partial_x^l Xu \big] \\
	=& (1-x^2) \partial_x^{l+2} \left( (1-x^2)^{l+1 -\frac{n}{2}} \partial_x^l Xu \right) \\
	& + (l+1) \cdot (-2x) \cdot \partial_x^{l+1} \left( (1-x^2)^{l+1 -\frac{n}{2}} \partial_x^l Xu \right) \\
	& + \frac{(l+1)l}{2} (-2) \cdot \partial_x^{l} \left( (1-x^2)^{l+1 -\frac{n}{2}} \partial_x^l Xu \right) \\
	& -(n-2l-2)x \cdot \partial_x^{l+1} \left( (1-x^2)^{l+1 -\frac{n}{2}} \partial_x^l Xu \right) \\
	& - (l+1)(n-2l-2) \cdot \partial_x^{l} \left( (1-x^2)^{l+1 -\frac{n}{2}} \partial_x^l Xu \right) \\
	=&  (1-x^2)\partial_x^2 e_l u  - nx\partial_x e_l u - (l+1) (n-l-2) e_l u \\
	= & (\Delta_{\Sn} - (l+1)(n-l-2)) e_l u\\
	=& e_{l+1} u.
\end{align*}
This finishes the induction.

Finally, we conclude the proof by
 $$P_n u =(-1)^{\frac{n}{2}} e_{\frac{n}{2} - 1} u = (-1)^{\frac{n}{2}} \partial_{n-2} X u = (-1)^{\frac{n}{2}}  \partial_x^{n-1} \left( (1-x^2)^{\frac{n}{2}} \partial_x u \right).$$ 
\end{proof}

\medskip

\noindent{\bf AI disclosure.} The authors knew that the nonlocal factor of $P_{2m+1}$ on $\Sp^{2m+1}$ is closely related to the Dirichlet-to-Neumann map before September 2026. We use AI to find its relevance to the Hilbert transform, with a proof of \eqref{oddGJMS} in the spirit of \cite{Chen-Shi} using Gegenbauer polynomials. Based on this proof, we give a unified induction proof adapted for all dimensions. Another induction proof of \eqref{evenGJMS}
in Appendix \ref{Append} is due to the second named author in 2025 without AI.

\bibliographystyle{amsplain}
\bibliography{references}

\end{document}